\documentclass[journal]{IEEEtran}

\usepackage{tikz}

\usepackage{amsmath}
\usepackage{latexsym, amssymb}
\usepackage{graphicx}
\usepackage{amsmath, amsbsy}
\usepackage{amsopn, amstext}
\usepackage{hyperref}
\usepackage{cancel, color}
\usepackage{epstopdf}
\usepackage{float}

\newcommand{\bea}{\begin{eqnarray}}
\newcommand{\eea}{\end{eqnarray}}
\def\nn{\nonumber}

\DeclareMathOperator{\Span}{Span}
\DeclareMathOperator{\Col}{Col}

\DeclareMathOperator{\lcm}{lcm}
\def\cal{\mathcal}

\def\diag{diag}
\def\ra{\rightarrow}

\def\d{\delta}

\def\D{\Delta}

\def\0{{\bf 0}}

\newcommand{\R}{{\mathbb R}}

\def\dsum{\mathop{\sum}\limits}

\newtheorem{thm}{Theorem}[section]

\newtheorem{cor}[thm]{Corollary}
\newtheorem{dfn}[thm]{Definition}
\newtheorem{prp}[thm]{Proposition}
\newtheorem{exa}[thm]{Example}
\newtheorem{rem}[thm]{Remark}

\begin{document}

\title{Cooperative Control via Congestion Game Approach}


\author{Yaqi~Hao,~Sisi~Pan,~Yupeng~Qiao,
	and~Daizhan~Cheng, ~\IEEEmembership{Fellow,~IEEE}
	\thanks{This work is supported partly by the National Natural Science Foundation of China (NSFC) under Grant 61333001.}
	\thanks{Yaqi Hao is with Control Science and Engineering, Shandong University, Ji'nan 250061, P.R.China e-mail: hoayaqi@outlook.com.}
    \thanks{Yupeng Qiao and Sisi Pan are with College of Automation Science and Engineering, South
China University of Technology, Guangzhou 510640, P.R.China e-mail: ypqiao@scut.edu.cn.}
    \thanks{Daizhan Cheng is with Key Laboratory of Systems and  Control, Chinese Academy of Sciences, Beijing 100190, P. R. China e-mail: dcheng@iss.ac.cn.}
	\thanks{ Corresponding author: Daizhan Cheng. Tel.: +86 1062651445; fax.: +86 1062587343.}
}

\markboth{IEEE Transactions On Automatic Control, Vol. XX, No. Y, Month 201Z}
{Hao \MakeLowercase{\textit{et al.}}: Congestion Game Approach}

\maketitle

\begin{abstract}
The optimization of facility-based systems is considered. First, the congestion game is converted into a matrix form, so that the matrix approach is applicable. Then, a facility-based system with a system performance criterion is considered. A necessary and sufficient condition is given to assure that the system is convertible into a congestion game with the given system performance criterion as its potential function by designing proper facility-cost functions. Using this technology, for a dynamic facility-based system the global optimization may be reached when each agent optimizes its payoff functions. Finally, the approach is extended to those systems which are partly or nearly convertible.
\end{abstract}


\IEEEpeerreviewmaketitle

\section{Introduction}

The distributed resource allocation problem, such as distributed welfare \cite{mar10}, cost sharing \cite{che09,fal12}, etc., aims at optimal resource distribution.
This problem has been formulated as a congestion game, which is a special class of potential games \cite{ros73,mon96b}. Precisely speaking, by designing proper utility functions to each agent, the overall welfare (or overall cost) is considered as the potential function. Then, the techniques developed for game theoretic control (GTC) are applicable to find pure Nash equilibriums, which provide candidates of optimal solutions \cite{gop11,gop14}.

In a distributed resource allocation problem, the overall welfare is separable, if it can be written as
\begin{align}\label{0.1}
W=\dsum_{r\in R}W_r,
\end{align}
where $\{W_r\}$ is the set of separated welfare functions. In recent works this separability is assumed.
For instance, as described in \cite{mar10}, distributed welfare game is a tuple $G=\{N,R,\{{\cal A}_i\}, \{W_r\},\{f_r\}\}$, where $\{W_r\}$, the welfare function for resource $r$, is known and the overall welfare $W=\dsum_{r\in R}W_r$ is determined.

In this note we consider a facility-based system as $G=\left\{N,R,\{{\cal A}_i\}, P\right\}$, where $P=W$ is a global performance criterion, which could be considered as overall welfare/cost, etc. The main problem concerned here is: can we convert this facility-based system to a congestion game? That is, can we design the facility-cost functions $\Xi_r$, such that the corresponding welfare functions $P_r=W_r$ for facility $r$ satisfy (\ref{0.1})? Briefly speaking, we want to know whether $W$ is separable? If ``yes", the GTC techniques developed in \cite{mar10,che09,fal12}, etc., can be used for facility-based systems.

In resorting to semi-tensor product (STP) of matrices, the problem is investigated by expressing a congestion game into its matrix form. Then,  the separability problem becomes solving a set of linear equations. The main contribution of this note consists of two parts: (1) Check whether the objective function is separable. If ``yes", the design of cost functions is proposed. (2) If ``no", the nearest separable potential game is considered, which enlarges the applicable set of  the previous design method to facility-based games.

 The STP of matrices is a generalization of conventional matrix product and all the computational properties of conventional matrix product remain available. Throughout this note, the default matrix product is STP, so the product of two arbitrary matrices is well defined and the symbol $\ltimes$ is mostly omitted. A brief survey on STP and related notations are provided in Appendix.

The rest of this note is organized as follows: Section 2 proposes a matrix form description for a congestion game. Section 3 provides a necessary and sufficient condition for the separability. In Section 4, the congestion game approach is extended to cases where facilities are either restricted or inconsistent.

\section{Matrix Expression of Congestion Games}

A congestion game is a tuple $G=\left(N,M,(A^i)_{i\in N}, (\Xi_j)_{j\in M}\right)$, where $N$ is the set of players; $M$ is the finite set of facilities to be shared by players; the facility-cost function $\Xi_j:~\R\rightarrow \R$ describes the cost of facility $j\in M$, which depends on the number of players using the facility $j$ in a profile $a$; $A^i\subset 2^M$ is the strategy (action) set of player $i$ and each strategy (action) in $A^i$ is a subset of $M$, which means that player $i$ has the option of selecting multiple facilities \cite{mon96b}.

Denote the set of profiles as $A=\prod_{i=1}^nA^i$. For a profile $a=(a^1,\cdots,a^n)\in A$ the number of users of facility $j$ is denoted as
\begin{align}\label{2.1}
r_j(a):=\left|\{i\;\big|\;j\in a^i\}\right|,\quad j=1,\cdots,m.
\end{align}
Define the payoff function of player $i$, i.e., $c_i:A\ra \R$ by
\begin{align}\label{2.2}
c_i(a):=\dsum_{j\in a^i}\Xi_j(r_j(a)),\quad i=1,\cdots,n,
\end{align}
and a function $P:A\ra \R$ as
\begin{align}\label{2.3}
P(a):=\dsum_{j\in \cup_{i=1}^na^i}\left(\dsum_{\ell=1}^{r_j(a)} \Xi_j(\ell)\right).
\end{align}
Using (\ref{2.2}) and (\ref{2.3}), we have the following result:

\begin{thm}\label{th2.2} \cite{mon96b} A congestion game is a potential game with payoff functions in (\ref{2.2}) and the potential function in (\ref{2.3}).
\end{thm}

In the following, we will express (\ref{2.2}) and (\ref{2.3}) into matrix forms.

Assume
$$
N=\{1,2,\cdots,n\};\quad M=\{1,2,\cdots,m\}.
$$
To begin with, we consider the facility-cost functions $\Xi_j$. Denote
$$
\Xi_j(k):=\xi^j_k,\quad k=1,\cdots,n;\;j=1,\cdots,m,
$$
where $k$ is the number of players using facility $j\in M.$
Then, $\Xi_j$ can be expressed in a vector form as
$$
\Xi_j=\left[\xi^j_1,\xi^j_2,\cdots,\xi^j_n\right],\quad j=1,\cdots,m.
$$
Putting all facility-cost functions together, we have
\begin{align}\label{2.4}
\Xi:=\left[\Xi_1,\Xi_2,\cdots,\Xi_m\right]\in \R^{mn}.
\end{align}

Next, we express each strategy (action)  $a^i\in A^i$ into a vector form. Since $a^i\in 2^M$, we use an index vector to express it. Let $a^i\in \R^m$ be a column vector with entries as
\begin{align}\label{2.5}
a^i(s):=\begin{cases}
1,\quad s\in a^i;\\
0,\quad \mbox{otherwise}.
\end{cases}
\end{align}
Using (\ref{2.1}), we construct $$r(a):=[r_1(a),\cdots,r_m(a)]^{\mathrm{ T } }.$$
It can be verified that
\begin{align}\label{2.6}
r(a)=\dsum_{i=1}^na^i,\quad a\in A.
\end{align}

Define
$$
d_i(a)=
\begin{cases}
\d_n^{r_i(a)},\quad r_i(a)\neq 0;\\
{\bf 0}_n,\quad \mbox{otherwise}.
\end{cases}
$$
Using them, we construct
\begin{align}\label{2.7}
D(a):=\diag\left(d_1(a),d_2(a),\cdots,d_m(a)\right).
\end{align}

A straightforward computation shows the following result:

\begin{prp}\label{p2.3}
The payoff functions $c_i$ can be expressed as
\begin{align}\label{2.8}
c_i(a)=\Xi D(a) a^i,\quad i=1,\cdots,n.
\end{align}
\end{prp}

Finally, we construct a set of Boolean vectors  as
$$
\begin{array}{l}
b_i(a):=\underbrace{[1,\cdots,1}_{r_i(a)},\underbrace{0,\cdots,0]}_{n-r_i(a)},\quad i=1,2,\cdots,m,
\end{array}
$$
and define
\begin{align}\label{2.9}
B(a):=\left[b_1(a),b_2(a),\cdots,b_m(a)\right].
\end{align}
Then, we can prove the following result.

\begin{prp}\label{p2.4} The potential function can be expressed as
\begin{align}\label{2.10}
P(a)=\Xi B^{\mathrm{T}}(a),\quad a\in A.
\end{align}
\end{prp}

\section{Optimization via Designed Facility Costs}

\subsection{Design of Facility-Cost Functions}

First, we give a rigorous description for a facility-based system (FBS):

\begin{dfn}\label{d3.1}
\begin{enumerate}
\item  An FBS is a tuple $\Sigma=(N, M, (A^i)_{i\in N},P)$. Here,  $N=\{1,2,\cdots,n\}$ is the set of players and $M=\{1,2,\cdots,m\}$ is the set of facilities shared by players. Each player is capable of selecting potentially multiple facilities in $M$; therefore, we say that   player $i$ has strategy (or action) set $A^i\subset 2^M$, that is, the set of certain subsets of $M$. $P:~A\ra \R$ is the system overall cost, which needs to be minimized.
\item $A:=\prod_{i=1}^nA^i$ is called the set of profiles of the system.
\end{enumerate}
\end{dfn}

Our purpose is to find a profile $a^*\in A$, such that
\begin{align}\label{3.1}
P(a^*)=\min_{a\in A}P(a).
\end{align}

The fundamental idea of the technique developed in this paper is: choosing suitable facility-cost functions (FCFs) such that the FBS can be converted into a congestion game with a pre-assigned performance criterion $P(a)$ as its potential function. Then, we use  properties of the potential game to realize the optimization.
Therefore, the key issue is: can we find a suitable set of FCFs such that the given $P(a)$ becomes its potential function? To answer this question, we construct a linear system as follows:

Assume $|A|=\ell$ (that is, there are $\ell$ different profiles) and denote $A=\{S_1,S_2,\cdots,S_{\ell}\}$. A linear system is defined as
\begin{align}\label{3.2}
B\Xi^{\mathrm{T}}=P,
\end{align}
where
$$
B=\begin{bmatrix}
B(S_1)\\
B(S_2)\\
\vdots\\
B(S_{\ell})
\end{bmatrix},\quad
P=\begin{bmatrix}
P(S_1)\\
P(S_2)\\
\vdots\\
P(S_{\ell})
\end{bmatrix},
$$
and $B(a)$ is defined in (\ref{2.9}).

The following result answers the above question.

\begin{thm}\label{t3.2} Consider an FBS. A set of FCFs can be found such that the FBS becomes a congestion game with $P(a)$ as its potential function, if and only if, (\ref{3.2}) has at least one solution.
\end{thm}

\begin{IEEEproof} The necessity comes from the matrix expression of a congestion game. Precisely speaking, collecting (\ref{2.10}) for all $a\in A$ together yields (\ref{3.2}). As for the sufficiency, using the solution of (\ref{3.2}), i.e., the $\Xi$, it is easy to construct the corresponding FCFs. Then a straightforward computation shows that the corresponding potential function is exactly the given $P(a)$.
\end{IEEEproof}

\begin{rem}\label{r3.3} Using (\ref{3.2}), we have
$$
\begin{array}{ccl}
P(a)&=&\Xi B^{\mathrm{T}}(a)\\
~&=&\Xi_1B_1^{\mathrm{T}}(a)+\cdots+ \Xi_mB_m^{\mathrm{T}}(a)\\
~&=&\left(\dsum_{j=1}^{r_1(a)}\Xi^1_j\right)+\cdots+\left(\dsum_{j=1}^{r_m(a)}\Xi^m_j\right)\\
~&=&\dsum_{i\in M}P_i(\#(a)_i),
\end{array}
$$
where
$$
P_i(\#(a)_i)=\dsum_{j=1}^{r_i(a)}\Xi^i_j
$$
is the separated cost for facility $i$, ~$i=1,\cdots,m$.

One can easily see that
$$
\Xi^i_k:=c_i(k)=P_i(k)-P_i(k-1).
$$
Hence, (\ref{3.2}) implies a standard distributed cost structure \cite{gop11}.\footnote{This is pointed out by an anonymous reviewer.}
\end{rem}

\subsection{Dynamics of Facility-Based Systems}

\begin{dfn}\label{d4.1} An FBS is called a dynamic FBS, if the system (or game) is repeated and the strategy profile is updated in a Markov-style as:
\begin{align}\label{3.4}
\begin{cases}
x_1(t+1)=f_1(x_1(t),x_2(t),\cdots,x_n(t))\\
x_2(t+1)=f_2(x_1(t),x_2(t),\cdots,x_n(t))\\
\quad\quad\quad\quad\ \vdots\\
x_n(t+1)=f_n(x_1(t),x_2(t),\cdots,x_n(t)).
\end{cases}
\end{align}
\end{dfn}

The dynamic equation (\ref{3.4}) is determined by the strategy updating rules (SURs). There are several commonly used SURs. We refer to \cite{che15}, \cite{Qi16} for more SURs and the constructing process of building the dynamic model using SURs. Only one SUR, called the Myopic Best Response Arrangement (MBRA), is used in this paper. We briefly describe it as follows: Consider a game with $N=\{1,2,\cdots,s\}$ and player $i$ has its strategies as $A^{i}=\{1,\cdots,k_i\}$. Assume at time $t$ the other players use their strategies as $a^{-i}\in \prod_{j\neq i}A^{j}$, and the player $i$ is allowed to update his strategy at the next moment $t+1$, then he will choose
\begin{align}\label{3.5}
x_i(t+1)\in \arg min_{\ell \in A^{i}}\left(c_i(\ell,a^{-i})\right).
\end{align}
MBRA is widely used because it has the following nice property.

\begin{thm}\label{t3.4}\cite{mon96b} Consider a finite potential game.
\begin{enumerate}
\item It has at least one pure Nash equilibrium.
\item If at each moment only one player is allowed to update his strategy and MBRA is used, then the dynamic potential game will converge to its Nash equilibrium.
\end{enumerate}
\end{thm}

\begin{rem}\label{r3.5}
\begin{enumerate}
\item Assume (\ref{3.2}) has solution. The profile $a^{\ast},$ such that (\ref{3.1}) holds, is a Nash equilibrium.
\item It is worth noting that in general a congestion game may have more than one Nash equilibriums. As long as the Nash equilibrium is unique, MBRA-based dynamics will converge to this unique Nash equilibrium, which is also the minimum point of $P(a)$.
\end{enumerate}
\end{rem}

\section{Design for Restricted FBSs}

\subsection{Partly Designable Facilities}
This subsection considers the case when only part of facility-cost functions  can be designed. This situation happens, for instance, the rest facilities are owned by other companies or so, and hence you can only design the facility-cost functions of your own facilities.

Based on the attribution of design right, the finite set of facilities $M$ can be divided into two disjoint sets as follows:
$$M=\Omega\cup\Omega_{\complement},$$
such that for each facility $j,$ its facility-cost function can be designed, if and only if, $j\in \Omega_{\complement}.$

Now set $\Omega=\{i_{1},i_{2},\cdots,i_{t}\}\subseteq M.$ According to this partition in $M$, similarly, we divide $\Xi,$ $B(a)$ and $P(a)$ respectively as
\begin{equation}
\begin{array}{ccl}
\tilde{\Xi}&:=&[\Xi_{1},\cdots,\Xi_{i_{1}-1},\Xi_{i_{1}+1},\cdots,\Xi_{i_{t}-1},\nn\\
&&\Xi_{i_{t}+1},\cdots,\Xi_m],\nn\\
\hat{\Xi}&:=&[\Xi_{i_{1}},\Xi_{i_{2}},\cdots,\Xi_{i_{t}}],\nn\\
\tilde{B}(a)&:=&[b_{1}(a),\cdots,b_{i_{1}-1}(a),b_{i_{1}+1}(a),\cdots,\nn\\
&&\ b_{i_{t}-1}(a),b_{i_{t}+1}(a),\cdots,b_{m}(a)],\nn\\
\hat{B}(a)&:=&[b_{i_{1}}(a),b_{i_{2}}(a),\cdots,b_{i_{t}}(a)],\nn\\
\tilde{P}(a)&:=&P(a)-\hat{B}(a)\hat{\Xi}^{\mathrm{T}}\nn.
\end{array}
\end{equation}
From above definitions, one can easily see that $\tilde{\Xi},~\tilde{B}(a)$ and $\tilde{P}(a)$ are only related to facilities in $\Omega_{\complement}.$ Using above denotations, a linear system is defined as

\begin{align}\label{4.2}
\tilde{B}\tilde{\Xi}^{\mathrm{T}} =\tilde{P},
\end{align}
where
$$
\tilde{B}=\begin{bmatrix}
\tilde{B}(S_1)\\
\tilde{B}(S_2)\\
\vdots\\
\tilde{B}(S_{\ell})
\end{bmatrix},\quad
\tilde{P}=\begin{bmatrix}
\tilde{P}(S_1)\\
\tilde{P}(S_2)\\
\vdots\\
\tilde{P}(S_{\ell})
\end{bmatrix}.
$$

The following result is an immediate consequence of Theorem  \ref{t3.2}.

\begin{cor}\label{c4.2} Consider an FBS with a given performance criterion $P(a)$. If only  part of facility-cost functions can be designed, a set of FCFs can be found such that the FBS becomes a congestion game with $P(a)$ as its potential function, if and only if, (\ref{4.2}) has at least one solution.
\end{cor}
\subsection{Restricted Facilities}
This subsection  considers the case when some facilities are restricted because of, for instance, their ultimate bearing capacity. Hence, you can only consider desirable profiles.

Assume a set of constraints are given as
\begin{align}\label{4.5}
\begin{cases}
k_{11}r_1(a)+k_{12}r_2(a)+\cdots+k_{1m}r_m(a)<T_{1}\\
k_{21}r_1(a)+k_{22}r_2(a)+\cdots+k_{2m}r_m(a)<T_{2}\\
\quad\quad\quad\quad\quad\quad\quad\quad\quad\vdots\\
k_{t1}r_1(a)+k_{t2}r_2(a)+\cdots+k_{tm}r_m(a)<T_{t}.
\end{cases}
\end{align}

We say, for a profile $a$, it is a desirable one, if and only if, it satisfies (\ref{4.5}). Otherwise, it is undesirable.

Hence, according to (\ref{4.5}), we can verify that the set of profiles $A$ contains two disjoint parts expressed as
$$A=\Omega\cup \Omega^{c},$$
where $\Omega$ is the set of desirable profiles while $\Omega^{c}$ is the set of undesirable ones.

In order to assure the Nash equilibrium $a^{\ast}\in \Omega,$ we further define payoff functions $\widetilde{c}_{i}(a)$ and performance criterion $\widetilde{P}(a)$ respectively as follows:
\begin{align}\label{4.9}
\tilde{c}_{i}(a):=
\begin{cases}
c_{i}(a),\quad a\in\Omega; \\
c^{\ast},\ \quad\quad a\in\Omega^{c},
\end{cases}
\end{align}
where $c^{\ast}\gg \max \limits_{i=1,2,\cdots,n}\{c_{i}(a)|a\in \Omega\};$

\begin{align}\label{4.12}
\tilde{P}(a):=
\begin{cases}
P(a),\quad a\in\Omega;  \\
P^{\ast},\quad\quad a\in\Omega^{c},
\end{cases}
\end{align}
where $P^{\ast}\gg \max\{P(a)|a\in \Omega\}$.

From the definition of $\Omega$ and (\ref{4.12}),  the following result can be obtained.

\begin{prp}\label{p4.1} The minimization of performance criterion $P(a)$ is equivalent to that of $\tilde{P}(a)$. That is,
$$
\begin{array}{ccl}
 P(a^{\ast})=\min_{a\in A}P(a)
&\Leftrightarrow& P(a^{\ast})=\min_{a\in \Omega}P(a)\\
~&\Leftrightarrow& \tilde{P}(a^{\ast})=\min_{a\in A}\tilde{P}(a).
\end{array}
$$
\end{prp}

A linear system is defined as
\begin{align}\label{4.3}
B\Xi^{\mathrm{T}}=P,
\end{align}
where$$
B=\begin{bmatrix}
B(S_{l_{1}})\\
B(S_{l_{2}})\\
\vdots\\
B(S_{l_{k}})
\end{bmatrix},\quad
P=\begin{bmatrix}
P(S_{l_{1}})\\
P(S_{l_{2}})\\
\vdots\\
P(S_{l_{k}})
\end{bmatrix},\quad\{S_{l_{1}},\cdots,S_{l_{k}}\}=\Omega.
$$
According to Proposition \ref{p4.1} and Theorem \ref{t3.2}, we have the following result.

\begin{cor}\label{c4.4} Consider an FBS with a given performance criterion $P(a)$. If facilities are restricted as (\ref{4.5}), a set of FCFs can be found such that the FBS becomes a congestion game with $P(a)$ as its potential function, if and only if, (\ref{4.3}) has at least one solution.
\end{cor}
\subsection{FBS with Improper $P$}

Consider an FBS $G$ with given $\Xi$ and $P$. Assume for this given $P$  (\ref{3.2}) has no solution. For this case, the dynamical equivalence is applied to investigate its convergence to Nash equilibriums.

\begin{dfn}\cite{che15.2}\label{d5.1} Two evolutionary games are said to be dynamically equivalent, if they have the same strategy profile dynamics (that is, $f_{i},~i=1,2,\cdots,n,$ defined in (\ref{3.4})).
\end{dfn}

First, we give the process of constructing an FBS $G's$ closest congestion game $G_{0}.$

Let $B_0$ be the matrix consisting of maximum linear independent columns of $B$ defined in (\ref{3.2}). Hence, $B_0$ has full column rank.

From (\ref{3.2}) one sees easily that the subspace of congestion games, denoted by $V^p$, is spanned by the columns of $B_0$. That is,
\begin{equation}
V^p=\Span\{Col(B_0)\}.
\end{equation}
For an FBS $G$ with given $P(a)$, it is clear that the least square solution for (\ref{3.2}) is
\begin{equation}\label{5.9}
  \Xi_0^{\mathrm{T}}  =(B_0^{\mathrm{T}} B_0)^{-1}B_0^{\mathrm{T}}P.
\end{equation}
Using this $\Xi_0$, the corresponding potential function $P_0$ can be calculated via (\ref{3.2}).
Then, we have the following definition.
\begin{dfn}
For an FBS $G$ with given $P$, the game $G_{0}$ is said to be its closest congestion game, if the facility-cost functions and potential function are
$\Xi_{0}$ and $P_{0}$ mentioned above.
\end{dfn}

According to Definition \ref{d5.1}, the following proposition is straightforward verifiable.
\begin{prp}\label{p5.2}
If an FBS $G$ with given $P$ and $\Xi$ is dynamically equivalent to its closest congestion game $G_{0}$, then it can be led to a pure Nash equilibrium.
\end{prp}

Assume there exists $\epsilon\geq 0,$ such that
\begin{align}\label{5.11}
\|P-P_0\|:=\max_{a\in A}|P(a)-P_0(a)|<\epsilon.
\end{align}
The following result can be obtained.

\begin{prp}\label{p5.4}  Assume (\ref{5.11}) and MBRA is used to generate the dynamics of $G$ by using $\Xi$ and $\Xi_0$ respectively. If these two dynamics are dynamically equivalent, then $G$ with $\Xi$ will converge to a Nash equilibrium. If the Nash equilibrium, denoted by $a^*$, is unique, it is a near smallest point, that is,
$$
|P(a^*)-P_{\min}| < 2\epsilon.
$$
\end{prp}

\begin{exa}\label{e5.5}

Given an FBS with $N=\{1,2,3\}$, $M=\{1,2,3,4,5\}$,
$$
\begin{array}{l}
A^1=\{a_{1}^{1},a_{2}^{1}\},\\
A^2=\{a_{1}^{2},a_{2}^{2},a_{3}^{2}\},\\
A^3=\{a_{1}^{3},a_{2}^{3},a_{3}^{3}\},
\end{array}
$$
where $a^1_1=\{1,2,3\},\ a^1_2=\{3,4,5\},$
      $a^2_1=\{1,2,4\},\ a^2_2=\{3,5\},\ a^2_3=\{4,5\},$
      $a^3_1=\{1,3,4\},\ a^3_2=\{2,5\},\ a^3_3=\{3,5\}.$

Denote the set of profiles (in alphabetic order) as
$\{S_1,S_2,\cdots,S_{18}\}=\{a^1_1a^2_1a^3_1,~a^1_1a^2_1a^3_2,~\cdots, ~a^1_2a^2_3a^3_3\}.$ Using (\ref{2.6}), we have
$$
\begin{array}{l}
r(S_1)=[3,2,2,2,0]^{T},\\
r(S_2)=[2,3,1,1,1]^{T},\\
\quad\quad\quad\vdots\\
r(S_{18})=[0,0,2,2,3]^{T}.\\
\end{array}
$$
It follows from (\ref{2.9}) that
$$
B(S_1)=[1,1,1,1,1,0,1,1,0,1,1,0,0,0,0].
$$
Similarly, $B(S_i)$, $i=2,3,\cdots,18,$ can also be calculated. Then the coefficient matrix $B$ for (\ref{3.2}) is obtained as
$$\setlength{\arraycolsep}{0.2pt}
B=
   \begin{pmatrix}
        \begin{smallmatrix}
    1 & 1 & 1 & 1 & 1 & 0 & 1 & 1 & 0 & 1 & 1 & 0 & 0 & 0 & 0\\
    1 & 1 & 0 & 1 & 1 & 1 & 1 & 0 & 0 & 1 & 0 & 0 & 1 & 0 & 0\\
    1 & 1 & 0 & 1 & 1 & 0 & 1 & 1 & 0 & 1 & 0 & 0 & 1 & 0 & 0\\
    1 & 1 & 0 & 1 & 0 & 0 & 1 & 1 & 1 & 1 & 0 & 0 & 1 & 0 & 0\\
    1 & 0 & 0 & 1 & 1 & 0 & 1 & 1 & 0 & 0 & 0 & 0 & 1 & 1 & 0\\
    1 & 0 & 0 & 1 & 0 & 0 & 1 & 1 & 1 & 0 & 0 & 0 & 1 & 1 & 0\\
    1 & 1 & 0 & 1 & 0 & 0 & 1 & 1 & 0 & 1 & 1 & 0 & 1 & 0 & 0\\
    1 & 0 & 0 & 1 & 1 & 0 & 1 & 0 & 0 & 1 & 0 & 0 & 1 & 1 & 0\\
    1 & 0 & 0 & 1 & 0 & 0 & 1 & 1 & 0 & 1 & 0 & 0 & 1 & 1 & 0\\
    1 & 1 & 0 & 1 & 0 & 0 & 1 & 1 & 0 & 1 & 1 & 1 & 1 & 0 & 0\\
    1 & 0 & 0 & 1 & 1 & 0 & 1 & 0 & 0 & 1 & 1 & 0 & 1 & 1 & 0\\
    1 & 0 & 0 & 1 & 0 & 0 & 1 & 1 & 0 & 1 & 1 & 0 & 1 & 1 & 0\\
    1 & 0 & 0 & 0 & 0 & 0 & 1 & 1 & 1 & 1 & 1 & 0 & 1 & 1 & 0\\
    0 & 0 & 0 & 1 & 0 & 0 & 1 & 1 & 0 & 1 & 0 & 0 & 1 & 1 & 1\\
    0 & 0 & 0 & 0 & 0 & 0 & 1 & 1 & 1 & 1 & 0 & 0 & 1 & 1 & 1\\
    1 & 0 & 0 & 0 & 0 & 0 & 1 & 1 & 0 & 1 & 1 & 1 & 1 & 1 & 0\\
    0 & 0 & 0 & 1 & 0 & 0 & 1 & 0 & 0 & 1 & 1 & 0 & 1 & 1 & 1\\
    0 & 0 & 0 & 0 & 0 & 0 & 1 & 1 & 0 & 1 & 1 & 0 & 1 & 1 & 1\\
  \end{smallmatrix}
    \end{pmatrix}
.$$

\begin{enumerate}
\item Assume the system performance criterion $P(a)$ is given in Table \ref{Tab9.1}.

\begin{table}[!htbp] 
\centering \caption{Performance Criterion $P(a)$:\label{Tab9.1}}
\begin{tabular}{|c|c|c|c|c|c|c|c|c|c|}
\hline
$a$&111&112&113&121&122&123&131&132&133\\
\hline
$P(a)$ &33&27&24&26&23&25&25&22&20\\
\hline\hline
$a$&211&212&213&221&222&223&231&232&233\\
\hline
$P(a)$ &28&28&26&33&13&20&29&16&19\\
\hline
\end{tabular}
\end{table}

It is easy to verify that (\ref{3.2}) has solutions. For instance, one of the solutions is
\begin{align}\label{3.3}
\Xi=[11,2,4,0,5,6,0,3,7,2,6,3,1,3,4].
\end{align}
According to Theorem \ref{t3.2}, under the facility-cost functions determined by (\ref{3.3}), the system becomes a congestion game with pre-assigned $P(a)$ as its potential function.

\item
Assume the system performance criterion $P(a)$ and facility-cost functions $\Xi$ are given as
$$\begin{array}{l}
P=[29,25,24,28,12,18,25,24,19,27,29,24,32,19,\nn\\
\quad\quad~27,25,23,22]^T,\nn\\
\Xi=[0.5,0,0.5,1.5,5,2,5,0.5,10,11,5,3,0,0.5,0].\nn
\end{array}$$
It is easy to verify that for the given P (\ref{3.2}) has no solution. Hence, we consider its closest congestion game.

First, we can easily obtain $B_0$ from $B$ by deleting the last three columns of $B$ as
$$\setlength{\arraycolsep}{0.2pt}
B_{0}=
   \begin{pmatrix}
        \begin{smallmatrix}
   1 & 1 & 1 & 1 & 1 & 0 & 1 & 1 & 0 & 1 & 1 & 0\\
   1 & 1 & 0 & 1 & 1 & 1 & 1 & 0 & 0 & 1 & 0 & 0\\
   1 & 1 & 0 & 1 & 1 & 0 & 1 & 1 & 0 & 1 & 0 & 0\\
   1 & 1 & 0 & 1 & 0 & 0 & 1 & 1 & 1 & 1 & 0 & 0\\
   1 & 0 & 0 & 1 & 1 & 0 & 1 & 1 & 0 & 0 & 0 & 0\\
   1 & 0 & 0 & 1 & 0 & 0 & 1 & 1 & 1 & 0 & 0 & 0\\
   1 & 1 & 0 & 1 & 0 & 0 & 1 & 1 & 0 & 1 & 1 & 0\\
   1 & 0 & 0 & 1 & 1 & 0 & 1 & 0 & 0 & 1 & 0 & 0\\
   1 & 0 & 0 & 1 & 0 & 0 & 1 & 1 & 0 & 1 & 0 & 0\\
   1 & 1 & 0 & 1 & 0 & 0 & 1 & 1 & 0 & 1 & 1 & 1\\
   1 & 0 & 0 & 1 & 1 & 0 & 1 & 0 & 0 & 1 & 1 & 0\\
   1 & 0 & 0 & 1 & 0 & 0 & 1 & 1 & 0 & 1 & 1 & 0\\
    1 & 0 & 0 & 0 & 0 & 0 & 1 & 1 & 1 & 1 & 1 & 0\\
    0 & 0 & 0 & 1 & 0 & 0 & 1 & 1 & 0 & 1 & 0 & 0\\
    0 & 0 & 0 & 0 & 0 & 0 & 1 & 1 & 1 & 1 & 0 & 0\\
    1 & 0 & 0 & 0 & 0 & 0 & 1 & 1 & 0 & 1 & 1 & 1\\
    0 & 0 & 0 & 1 & 0 & 0 & 1 & 0 & 0 & 1 & 1 & 0\\
    0 & 0 & 0 & 0 & 0 & 0 & 1 & 1 & 0 & 1 & 1 & 0\\
     \end{smallmatrix}
    \end{pmatrix}
.$$
Using (\ref{5.9}), we have
\begin{align}
\Xi_{0}=&[0.4704,0.1516,0.0004,1.5766,4.6840,1.1375,\nn\\
&5.7214,0.1263,9.5267,11.2585,5.0109,2.5485,\nn\\
&0,0,0],\nn
\end{align}
and through (\ref{3.2}), the corresponding potential function $P_0$ can be calculated as
$$\begin{array}{l}
P_0=[29,25,23.9887,28.8315,12.5786,17.4214,\nn\\
\quad\quad~~ 24.3156,23.7109,19.1532,26.8641,28.7218,\nn\\
\quad\quad~~ 24.1641,32.1142,18.6828,26.6329,25.1359,\nn\\
\quad\quad~~ 23.5674,22.1170]^T.\nn
\end{array}$$

Next, we can calculate the payoff functions of $G$ by using $\Xi$ and $\Xi_0$ respectively, which are shown in Table \ref{Tab9.7} and Table \ref{Tab9.8} (roundoff to 0.01).

Using the MBRA, we can  get the best responding strategies respectively. It is obvious that they have the same strategy updating dynamics $f_1$,  $f_2$ and $f_3$, which are shown in Table \ref{Tab9.9}.

According to Definition \ref{d5.1}, we know that these two dynamics are dynamically equivalent. Hence, the FBS with given $\Xi$ and $P$ has at least one Nash equilibrium.

Moreover, for the FBS\ $G$ with $\Xi$, we have a switched system as
\begin{align}\label{5.10}
x(t+1)=L_{\sigma} x(t),
\end{align}
where $\sigma\in \{1,2,3\}$,
$$\begin{array}{l}
L_1=\d_{18}[10,2,3,4,5,6,7,17,9,10,2,3,4,5,6,7,17,9],\nn\\
L_2=\d_{18}[7,5,6,7,5,6,7,5,6,16,14,18,16,14,18,16,\nn\\
\quad\quad\quad~~~14,18],\nn\\
L_3=\d_{18}[3,3,3,5,5,5,9,9,9,12,12,12,14,14,14,18,\nn\\
\quad\quad\quad~~~18,18].\nn
\end{array}$$
Now, we assume the probability $P(\sigma=i)=1/3$, $i=1,2,3$, and three initial profiles are randomly chosen for a Matlab simulation. The results are shown in Fig.\ref{Fig.4}.

From Fig.\ref{Fig.4}, one sees that the strategy profile dynamics, starting from any initial profile, will converge to the unique Nash equivibrium $\d_{18}^{5}\sim \{1,2,2\}$.

Finally, we assume $\epsilon=0.9.$ It is easy to calculate that
\begin{align}
\|P-P_0\|:=\max_{a\in A}|P(a)-P_0(a)|=0.8315<\epsilon.
\end{align}
According to Proposition \ref{p5.4}, we know that the Nash equilibrium $\d_{18}^{5}\sim \{1,2,2\}$ is a near smallest point, that is,
$$
|P(a^*)-P_{\min}|<2\epsilon.
$$
In fact, it can be verified that for this example we have $
|P(a^*)-P_{\min}|=0$.

\begin{table}[!htbp] 
\centering \caption{Payoff matric of G with $\Xi$:\label{Tab9.7}}
\begin{tabular}{|c|c|c|c|c|c|c|c|c|c|}
\hline
$ c\backslash a$ & 111 & 112 & 113 & 121 & 122 & 123 & 131 & 132 & 133 \\
\hline
$c_1$ & 6 & 7 & 5.5 & 11.5 & 6 & 12 & 2 & 10.5 & 2.5 \\
\hline
$c_2$ & 10.5 & 13 & 16 & 10 & 1 & 10.5 & 5 & 11.5 & 11.5 \\
\hline
$c_3$ & 6 & 2 & 0.5 & 21 & 5.5 & 10.5 & 5.5 & 5.5 & 1 \\
\hline\hline
$c\backslash a$ & 211 & 212 & 213 & 221 & 222 & 223 & 231 & 232 & 233 \\
\hline
$c_1$ & 3.5 & 10.5 & 6 & 15.5 & 11.5 & 21 & 4 & 10 & 5.5 \\
\hline
$c_2$ & 4.5 & 10.5 & 7 & 10.5 & 0.5 & 10 & 3.5 & 5 & 5 \\
\hline
$c_3$ & 3.5 & 5.5 & 1 & 15.5 & 1.5 & 10 & 4 & 1.5 & 0.5  \\
\hline
\end{tabular}
\end{table}
\begin{table}[!htbp] 
\centering \caption{Payoff matric of G with $\Xi_{0}$:\label{Tab9.8}}
\begin{tabular}{|p{0.48cm}|p{0.48cm}|p{0.48cm}|p{0.48cm}|p{0.48cm}|p{0.48cm}|p{0.48cm}|p{0.48cm}|p{0.48cm}|p{0.48cm}|}
\hline
$ c\backslash a$ & 111 & 112 & 113 & 121 & 122 & 123 & 131 & 132 & 133 \\
\hline
$c_1$ & 4.81 & 7.01 & 4.96 & 11.25 & 5.28 & 11.57 & 1.85 & 10.88 & 2.17 \\
\hline
$c_2$ & 9.70 & 12.55 & 16.09 & 9.53 & 0.13 & 9.53 & 5.01 & 11.26 & 11.26 \\
\hline
$c_3$ & 5.14 & 1.14 & 0.13 & 20.94 & 4.68 & 9.53 & 5.29 & 4.68 & 0.13 \\
\hline\hline
$c\backslash a$ & 211 & 212 & 213 & 221 & 222 & 223 & 231 & 232 & 233 \\
\hline
$c_1$ & 2.67 & 10.73 & 5.14 & 14.54 & 11.38 & 20.79 & 2.67 & 10.73 & 5.14 \\
\hline
$c_2$ & 4.28 & 10.17 & 7.06 & 9.53 & 0.13 & 9.53 & 2.55 & 5.01 & 5.01 \\
\hline
$c_3$ & 2.83 & 4.68 & 0.13 & 15.01 & 1.58 & 9.53 & 3.15 & 1.58 & 0.13  \\
\hline
\end{tabular}
\end{table}
\begin{table}[!htbp] 
\centering \caption{strategy updating dynamics:\label{Tab9.9}}
\begin{tabular}{|c|c|c|c|c|c|c|c|c|c|}
\hline
$ f\backslash a$ & 111 & 112 & 113 & 121 & 122 & 123 & 131 & 132 & 133 \\
\hline
$f_1$ & 2 & 1 & 1 & 1 & 1 & 1 & 1 & 2 & 1\\
\hline
$f_2$ & 3 & 2 & 2 & 3 & 2 & 2 & 3 & 2 & 2\\
\hline
$f_3$ & 3 & 3 & 3 & 2 & 2& 2 & 3 & 3 & 3\\
\hline\hline
$f\backslash a$ & 211 & 212 & 213 & 221 & 222 & 223 & 231 & 232 & 233 \\
\hline
$f_1$ & 2 & 1 & 1 & 1 & 1 & 1 & 1 & 2 & 1\\
\hline
$f_2$ & 3 & 2 & 3 & 3 & 2 & 3 & 3 & 2 & 3\\
\hline
$f_3$ & 3 & 3 & 3 & 2 & 2 & 2 & 3 & 3 & 3\\
\hline
\end{tabular}
\end{table}
\begin{figure}[H]
\centering
  \includegraphics[width=3.5 in]{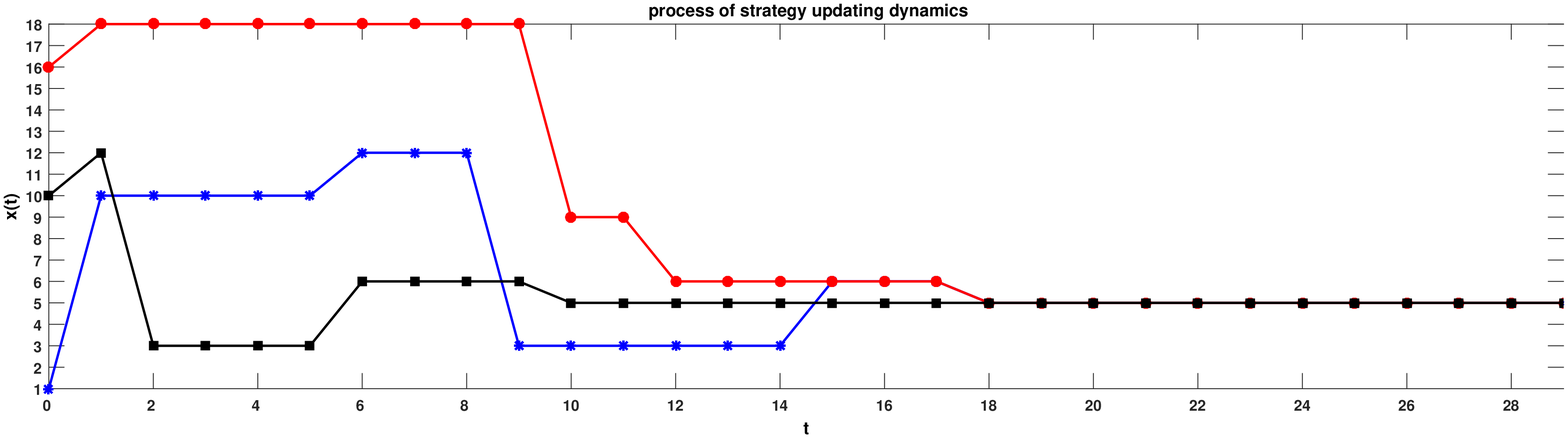}\\
  \caption{Profile Dynamics of (\ref{5.10})}\label{Fig.4}
\end{figure}
\end{enumerate}
\end{exa}

\section{Conclusion}

This note investigates the cooperative control of a facility-based system via congestion game approach. First, a matrix form description of a congestion game is presented. Then, a necessary and sufficient condition is obtained to assure the separability of the pre-assigned performance criterion $P(a)$. That is,  the FBS can be converted into a congestion game with $P(a)$ as its potential function. Using properties of potential games, the convergence to a Nash equilibrium is obtained. Thirdly, the result has been extended to some incomplete cases where $P(a)$ is not separable. Particularly, the problem of near dynamic congestion game is considered. It is proved that under dynamic equivalence the near congestion games may also be led to a Nash equilibrium.

Our approach can only be used for classical congestion games, where the user with multiple unit demands is not allowed. For instance, in the transportation congestion model a route segment can not be used by a player for more than once. But the user with multiple unit demands is an interesting and challenging problem. It could be studied in the future.

\section{Appendix}

This appendix gives a brief survey for semi-tensor product of matrices. We refer to \cite{che11,che12} for more details.

For technical statement ease, we first introduce some notations:
(1) ${\cal M}_{m\times n}$: the set of $m\times n$ real matrices.
(2) $\Col_i(M)$: the $i$-th column of $M$.
(3) ${\cal D}_k:=\left\{1,2,\cdots,k\right\},\quad k\geq 2$.
(4) $\d_n^i$: the $i$-th column of the identity matrix $I_n$.
(5) $\D_n:=\left\{\d_n^i\vert i=1,\cdots,n\right\}$.
(6) A matrix $L\in {\cal M}_{m\times n}$ is called a logical matrix
if the columns of $L$ are of the form of
$\d_m^k$.
Denote by ${\cal L}_{m\times n}$ the set of $m\times n$ logical
matrices.
(7) If $L\in {\cal L}_{n\times r}$, by definition it can be expressed as
$L=[\d_n^{i_1},\d_n^{i_2},\cdots,\d_n^{i_r}]$. For the sake of
compactness, it is briefly denoted as $
L=\d_n[i_1,i_2,\cdots,i_r]$.

The semi-tensor product of matrices is defined as follows:

\begin{dfn}\label{d1.1}  Let $M\in {\cal M}_{m\times n}$, $N\in {\cal M}_{p\times q}$, and $t=\lcm\{n,p\}$ be the least common multiple of $n$ and $p$.
The semi-tensor product (STP) of $M$ and $N$ is defined as
\begin{align}\label{1.1}
M\ltimes N:= \left(M\otimes I_{t/n}\right)\left(N\otimes I_{t/p}\right)\in {\cal M}_{mt/n\times qt/p},
\end{align}
where $\otimes$ is the Kronecker product.
\end{dfn}


In the following, we consider how to express a logical function into an algebraic form.

 Let $x_i\in {\cal D}_{k_i}$, $i=1,\cdots,n$ be logical variables and $f:\prod_{i=1}^n{\cal D}_{k_i}\ra {\cal D}_{k_0}$ be a (multi-valued) logical function. For ${\cal D}_k$, we identify $s\in {\cal D}_k$ with $\d_k^s\in \D_k$. Then we can express $f(x)$ into a matrix form.

 \begin{thm}\label{th1.3} \cite{che12} Let $x_i\in {\cal D}_{k_i}$, $i=1,\cdots,n$ be logical variables and $f:\prod_{i=1}^n{\cal D}_{k_i}\ra {\cal D}_{k_0}$ be a (multi-valued) logical function. When $x_i$ are expressed into vector form, there exists a unique logical matrix $M_f\in {\cal L}_{k_0\times k}$, where $k=\prod_{i=1}^nk_i$, such that
 \begin{align}\label{1.3}
 f(x_1,\cdots,x_n)=M_f\ltimes_{i=1}^nx_i.
 \end{align}
 $M_f$ is called the structure matrix of $f$.
\end{thm}

Next, assume a logical dynamic system
\begin{align}\label{1.4}
x_i(t+1)=f_i(x_1,\cdots,x_n),\quad i=1,\cdots,n.
\end{align}

Using Theorem \ref{th1.3}, (\ref{1.4}) can be expressed into matrix forms as
\begin{align}\label{1.5}
x_i(t+1)=M_ix(t),\quad i=1,\cdots,n,
\end{align}
where $M_i$ is the structure matrix of $f_i$, $i=1,\cdots,n$, and $x(t)=\ltimes_{i=1}^nx_i(t)$.
Multiplying both sides of (\ref{1.5}) together, we can have more compact form for (\ref{1.4}).

\begin{thm}\label{th1.4} \cite{che12} Consider system (\ref{1.4}). Using (\ref{1.5}), it can be expressed into its algebraic state space (ASS) form as
\begin{align}\label{1.6}
x(t+1)=Lx(t),
\end{align}
where $L=M_1*M_2*\cdots*M_n,$ and $*$ is the Khatri-Rao product.
\end{thm}

\end{document}